\documentclass[3p,onecolumn,11pt,sort&compress]{elsarticle}

\usepackage{amssymb}
\usepackage{amsmath}

\usepackage{float}

\usepackage{geometry}
\usepackage{setspace,caption}
\usepackage{chngcntr} 

\usepackage{enumerate}

\usepackage[symbol]{footmisc}
\usepackage{amssymb} 	
\usepackage[utf8]{inputenc}
\usepackage[T1]{fontenc}  
\usepackage{lmodern,textcomp}
\usepackage{amsmath}
\usepackage{mathtools}
\usepackage{resizegather}
\usepackage[yyyymmdd,hhmmss]{datetime}
\usepackage{amsfonts}
\usepackage{graphicx}
\usepackage[table]{xcolor}
\usepackage{physics}
\usepackage{gensymb}
\usepackage{array}
\usepackage{siunitx}
\usepackage[colorlinks=true, linkcolor=black, citecolor=blue, urlcolor=blue]{hyperref}
\usepackage{threeparttable}
\usepackage{booktabs}
\usepackage{subcaption}
\usepackage{multicol}
\usepackage{hhline}
\usepackage{diagbox}
\usepackage[absolute]{textpos}
\usepackage{multirow}
\usepackage{placeins}          
\usepackage[version=3]{mhchem} 
\usepackage[figuresright]{rotating}
\usepackage{lineno}
\usepackage{soul}
\usepackage{enumitem}
\usepackage{ltablex}
\usepackage{nicefrac}
\usepackage{cleveref}
\usepackage{braket}
\usepackage{lineno}

\crefname{chart}{Chart}{Charts}
\crefname{section}{Section}{Sections}
\crefname{figure}{Fig.}{Figs.}
\crefname{graph}{Graph}{Graphs}
\crefname{scheme}{Scheme}{Schemes}
\crefname{equation}{Eq.}{Eqs.}
\crefname{table}{Table}{Tables}
\crefname{chapter}{Chapter}{Chapters}
\crefname{appendix}{}{}

\journal{Applied Energy}

\begin{document}

\begin{frontmatter}

\title{Accurately modeling long-term storage with minimum representative hours in large-scale renewable energy systems} 

\author[rre]{Jacob Mannhardt}
\author[rre]{Lukas Kunz}
\author[rre]{Giovanni Sansavini\corref{corres}}
\cortext[corres]{Corresponding author: sansavig@ethz.ch}

\affiliation[rre]{organization={Institute of Energy and Process Engineering, ETH Zurich},
            city={Zurich},
            postcode={8092 Zurich}, 
            country={Switzerland}}

\begin{abstract}
Energy system optimization often relies on time series aggregation to ensure computational tractability. Aggregation generally loses the chronology of time steps, which renders the storage level representation challenging. Typically, this challenge is addressed by using representative days (RD) to utilize intra-day chronology, even though representative hours (RH) can describe the input time series more accurately at fewer representative time steps than RD. 
However, until now, the use of RH storage representation methods has been limited by either high computational complexity, poor accuracy in clustering and storage representation, or restricted applicability.
Here, we present a novel storage representation method based on RH that combines the high accuracy of RH time series aggregation with the high computational efficiency of methods based on RD. Through benchmarking the four most established storage representation methods on a model of a net-zero European energy system, we find that the proposed method can reduce the solving time by over 95\% for the same objective value compared to the most established RD and RH methods.
The proposed method exhibits particular strengths at strong aggregations of around 100 to 500 representative hours per year, making the method especially applicable to large-scale and sector-coupled transition pathway models.
The developed method for accurately modeling both short-term and long-term storage, along with the presented findings, is of practical relevance to energy system modelers who seek computational tractability in large-scale applications while avoiding the misallocation of storage and conversion capacities.
\end{abstract}




\end{frontmatter}

\section{Introduction}
Large-scale energy system optimization models (ESOMs) are key tools to outline the transition to low-carbon energy systems \cite{Victoria2022,Trondle2020a}. ESOMs face a conflict between accuracy, achieved by high spatiotemporal resolution and scope combined with technological detail, on the one hand, and computational tractability on the other. The literature deploys several techniques to improve the trade-off between accuracy and computational effort; namely, by aggregating the time series \cite{Gabrielli2018,Kotzur2018},
aggregating regions \cite{Frysztacki2023,Reinert2024}, selecting specific years on the transition pathway \cite{Mannhardt2024a,Victoria2022,Loffler2019} or optimizing a single year \cite{Trondle2020a,Neumann2023}, and simplifying technology characteristics \cite{Gabrielli2018a}. 

Time series aggregation (TSA) is the most established technique for complexity reduction \cite{Hoffmann2020,Nahmmacher2016,Teichgraeber2022}. TSA aggregates the data values of a fully resolved time series to a smaller number of representative clusters, thus reducing the size of the optimization problem. This aggregation can preserve the quantitative properties of the input time series, e.g., the total sum of the time series values \cite{Gonzato2021,Hoffmann2020}, but generally does not preserve the chronology of the representative time steps. As a result, modeling storage technologies poses a challenge because calculating the storage level involves coupling adjacent time steps. The challenge of modeling long-term storage is highlighted by the fact that many studies only include daily cycling without medium- to long-term storage, e.g., \cite{Scott2019, Backe2022,Buchholz2020,Fazlollahi2014}.
However, energy systems reliant on variable renewable resources require short- and long-term storage  \cite{Upadhyay2025,Guerra2020}. Omitting the accurate modeling of storage technologies can lead to strong misrepresentations of storage and conversion capacities. 

Hence, the storage representation in time-aggregated ESOMs must accurately capture both short-term and long-term storage while simultaneously supporting the reduced computational complexity of the TSA. To this end, the chronology of storage level time steps must be reconstructed to enable storage level coupling.
State-of-the-art storage representation methods \cite{Kotzur2018,Gabrielli2018,Pineda2018,Wogrin2016} allow for the inclusion of long-term storage, but always at the expense of computational complexity, accuracy, and applicability.

The two most common approaches for aggregating the hourly resolved time series to representative time steps are selecting representative days (RD) and selecting representative hours (RH) \cite{Hoffmann2021}. RD emerged because many energy system time series follow daily profiles, such as electricity and heat demand, and capacity factors of solar photovoltaic. Representative weeks have been applied at times \cite{Hoffmann2021}, but show limited effectiveness.
RD preserve intra-daily chronology, whereas RH are temporally independent of one another. Since RH are not constrained to form an intra-daily chronology, the original time series can be represented more accurately with fewer representative time steps, unless the time series only contain daily patterns. Therefore, RH are better suited to capture seasonal patterns with few representative time steps.

RD is the most established technique for TSA in ESOMs, \cite{Kotzur2018,Nahmmacher2016,Gabrielli2018,Jacobson2024,Novo2022,Parolin2024,Mantegna2024,VanDerHeijde2019,Martelli2020,Patil2021}, most of which deploy the superposition method by \citet{Kotzur2018} to model the storage level. RH has been utilized by fully resolving the storage level \cite{Gabrielli2018}, by using chronological RH (CRH) \cite{Pineda2018}, or by describing the transition between the system states \cite{Wogrin2016}. However, all three methods with RH show shortcomings in computational tractability \cite{Gabrielli2018}, representation of daily profiles \cite{Pineda2018}, approximation of the storage level \cite{Wogrin2016}, and applicability due to neglecting self-discharge \cite{Wogrin2016} (\cref{sec:lit_review}). 

In this paper, we introduce a novel storage representation method based on RH to utilize the accurate aggregation of the input time series when clustering RH, in comparison to RD and chronological RH. After clustering non-chronological RH, the chronology of the storage level time steps is reconstructed with as few storage level variables and constraints as possible to reduce computational complexity.
The presented approach significantly reduces the number of required variables and constraints without compromising any information, compared to the fully resolved storage level method by \citet{Gabrielli2018}. Therefore, this approach combines the accuracy of clustering RH with the computational efficiency of existing storage representation methods that use RD. Until now, the efficient representation of short- and long-term storage with self-discharge has not been achieved when using RH. 

We benchmark the presented approach against the original and another improved formulation by \citet{Kotzur2018} for RD, the approach by \citet{Gabrielli2018} for RH, and the approach by \citet{Pineda2018} for CRH. In this work, we focus solely on storage representation methods that enable the exact reconstruction of both short-term and long-term storage levels, thereby excluding the storage level approximation proposed by \citet{Wogrin2016}, which shares similarities with the proposed method but approximates the storage level evolution.
We find that the presented approach achieves high accuracy at low computational effort for a case study of a fully renewable European energy system model dependent on short- and long-term storage. In detail, we investigate the accuracy in objective value, conversion and storage technologies, and storage cycles.
We focus our analysis on strongly aggregated time series in large-scale ESOMs, i.e., reducing the number of representative time steps by up to 95\%, which is particularly relevant for pathway analyses. The presented storage representation method shows an objective error of only 7.5\% when heavily aggregating the original time series to 96 representative hours per year, while reducing the average solving time by more than 99\%.


The remainder of the paper is structured as follows: In \cref{sec:lit_review}, we provide a summary of existing methods forstorage representation under RD and RH. 
In \cref{sec:methods}, we briefly distinguish the aggregation to RD and RH (\cref{sec:clustering}), present the fully resolved storage constraints (\cref{sec:fully_resolved_storage_constraints}), and introduce the reduced storage representation with RH, presented in this paper (\cref{sec:met_ZEN_garden}). In \cref{sec:case_study}, the setup of the investigated case study is presented. In \cref{sec:results}, we benchmark the various storage representation methods in terms of objective value, solving time, and installed capacities. Finally, we assess the robustness of the investigated methods by applying them to three variations of the case study.
In \cref{sec:disc}, the findings are discussed and concluded.

\section{Storage representation in time-aggregated ESOMs - A short review}
\label{sec:lit_review}
The goal of time series aggregation (TSA) is to represent the set of input time series, i.e., of energy demands or renewable capacity factors, with a reduced number of time steps to reduce computational time. The solution complexity of state-of-the-art active-set (simplex) and interior-point (barrier) solvers generally increases polynomially or exponentially with the problem size, i.e., the number of variables and constraints \cite{Spielman2003,Wright1997}. Therefore, reducing the time resolution and, thereby, the number of variables and constraints can strongly decrease the solution time and memory requirements \cite{Hoffmann2021}. For large ESOMs with many regions, sectors, and a transition pathway of multiple years, the resulting optimization problem is typically computationally intractable otherwise. Common dimensionalities include hundreds of nodes \cite{Ganter2024}, up to 100 technologies \cite{Victoria2022}, and approximately 15 years in the transition pathway \cite{Mannhardt2024a}, albeit not combined in the same optimization problem.

The method for TSA predetermines the optimal method for modeling long-term storage (\cref{tab:lit_review}), since the relationships between representative time steps depend on the TSA method. 
In this work, we differentiate between \textit{representative days (RD)} and \textit{representative hours (RH)} as the entities used in the clustering for the identification of representative periods are days or hours, respectively \cite{Hoffmann2020}. For both approaches, the smallest increment, for which each time-dependent constraint and variable is formulated, is called \textit{representative time step}.
An RD contains 24 consecutive hours, and the RD are repeated in a non-chronological order. In contrast, using RH and selecting the representative hours instead of days generally does not preserve any chronology at all. A special case of RH are chronological RH (CRH), where a representative hour represents adjacent hours of the original time series, and thus the resulting RH are chronological. 

The most common methods to select the RD or RH are (i) clustering the time series, e.g., k-means or hierarchical clustering \cite{Nahmmacher2016,Hoffmann2020, Bahl2018, Teichgraeber2022}, and (ii) ``manual'' selection and downsampling, for instance, using every second hour or manually selected days per season \cite{Loulou2016,DeGuibert2020,Wiese2018,Gerbaulet2017}. For more details, see the extensive literature on clustering methods, e.g., Hoffmann et al. \cite{Hoffmann2020} and Teichgräber and Brandt \cite{Teichgraeber2022}.

RH have an additional degree of freedom over RD since they do not force the use of 24 chronological hours. As a result, clustering with RH represents the input time series more, or at worst equally, accurately for the same number of representative time steps \cite{Gonzato2021,Hoffmann2021}. 
Accordingly, RH can capture the properties of the full time series with fewer representative time steps. Since daily chronology is not enforced, RH can capture especially seasonal patterns better, e.g., heat demand or season-dependent renewable availability. On the other hand, RD need at least 48 representative time steps, i.e., two different representative days, to reflect the most basic seasonality. 

\begin{table}[htbp]
  \centering
  \caption{Commonly used time series aggregation and storage level representation approaches. Note that all storage level representation methods for representative hours (RH) could also be used with representative days (RD); however, the representation methods for RD cannot be sensibly applied to RH. The approach using fully resolved storage level variables and constraints by \citet{Gabrielli2018} is applied to both RD and RH. The original and additional references (ref.) are provided.}
\begin{tabular}{l|l}
\textbf{Time series aggregation} & \textbf{Storage level representation} \\
\midrule
\textit{Representative Days (RD)} & \textit{Superposition, hourly storage level limit constraint} \\
      & \multirow{2}[0]{*}{Original reference: \citet{Kotzur2018}} \\
Clusters of representative, &  \\
non-chronological days & Additional refs: \cite{Goke2022,Teichgraeber2022,Mantegna2024,Blanke2022} \\
\cmidrule{2-2}with 24 chronological hours each & \textit{Superposition, reduced storage level limit constraint} \\
      & \multirow{2}[0]{*}{Original reference: \citet{Kotzur2018}} \\
      &  \\
      & Additional refs: \cite{Gonzato2021,Parolin2024,Jacobson2024} \\
\cmidrule{2-2}      & \textit{Fully resolved storage level variable and constraint} \\
      & \multirow{2}[2]{*}{Original reference: \citet{Gabrielli2018}} \\
\cmidrule{1-1}\textit{Representative Hours (RH)} &  \\
      & Additional refs: \cite{Gonzato2021,Teichgraeber2022,Blanke2022,Raventos2020} \\
\cmidrule{2-2}Clusters of representative, & \textit{State-space formulation with transition matrix} \\
non-chronological hours & \multirow{2}[0]{*}{Original reference: \citet{Wogrin2016}} \\
      &  \\
      & Additional refs: \cite{Tejada-Arango2018,Tejada-Arango2018a} \\
\midrule
\textit{Chronological RH (CRH)} & \textit{Chronological storage level formulation} \\
      & \multirow{2}[0]{*}{Original reference: \citet{Pineda2018}} \\
Clusters of representative, &  \\
chronological hours & Additional refs: \cite{Goke2022,Teichgraeber2022,Gonzato2021} \\
\bottomrule
\end{tabular}%
  \label{tab:lit_review}%
\end{table}%

The literature proposes several methods to model long-term storage with RD (\cref{tab:lit_review}). The method initially presented by \citet{Gabrielli2018} fully resolves the storage level. In this case, the storage level is formulated for every hour of the full time series, and each hour of the full time series is uniquely mapped to an hour in the representative days. 

The main drawback of fully resolving the storage level is the number of required storage variables and constraints. As a result, the method has been found to even require more solving time than modeling the full time series \cite{Gonzato2021}. Therefore, most models follow the second method, introduced by \citet{Kotzur2018}. This method superposes an inter-daily storage level (between RD) with an intra-daily storage level (within RD).
Accordingly, the superposition reduces the number of required variables and constraints compared to the method by \citet{Gabrielli2018}. 
\citet{Kotzur2018} propose an improvement of their formulation, by constraining the storage level only twice per day to reduce the number of storage constraints, similarly to \cite{Parolin2024} and \cite{Gonzato2021}. 
The superposition of inter- and intra-daily levels is only sensibly applicable to RD clustering due to its use of chronological intra-daily storage levels. 
Hoffmann et al. \cite{Hoffmann2021} compare RD and RH clustering; however, they use the storage representation method of \citet{Kotzur2018} for both RD and RH in their analysis, which biases the results towards RD.

Although RH can represent the input time series more accurately, only a few works investigate storage modeling with RH (\cref{tab:lit_review}). The previously presented approach by \citet{Gabrielli2018} can use RD and RH, but the authors strictly use RD in their publication. \citet{Pineda2018} introduce a chronological aggregation and storage representation method, and compare it against representative days or weeks with daily or weekly cycling. Since the order of representative time steps is preserved in CRH, the storage level can be directly formulated with the same chronological time index. Chronological clustering and storage representation can strongly reduce the number of variables and constraints; however, the adjacency of represented hours smoothens daily and even weekly behavior. Intuitively, for less than 365 RH, no daily variation can be captured.
\citet{Gonzato2021} show that the chronological storage representation is not appropriate for a low number of RH. 
\citet{Wogrin2016} introduce a state-space model, where the sequence between RHs (or ``states'') is preserved by tracking the transitions between the RH. The storage level is not explicitly calculated, but it is implicitly tracked by summing up the changes in the storage level over time. Therefore, the state-space model cannot include self-discharge of the storage, as self-discharge depends on the absolute storage level \cite{Hall2008}. Furthermore, the transition between states is modeled as central differences, which approximates the storage level, and thereby does not ensure the correct tracking of the storage level.

\section{Methods}
\label{sec:methods}
\cref{sec:clustering} introduces the clustering of input time series into RD and RH. \cref{sec:fully_resolved_storage_constraints} introduces the constraints of storage technologies for fully resolved storage level time series. \cref{sec:met_ZEN_garden} presents the novel reduced storage representation method with RH introduced in this paper. 
For the sake of completeness, the other storage representation methods from the literature, which are used for benchmarking the presented approach, are detailed in the \cref{sec:storage_representation_RD,sec:storage_representation_RH}, using the same nomenclature as the presented method.

\begin{table}[htbp]
  \centering
  \caption{Summary of time step sets and mappings.}
\begin{tabular}{llll}
Description & Set & Index & Dimensionality  \\
\midrule
Full set of time steps & $\mathcal{T}$ & $t$ & $T=8760$ \\
Hours per day & $\mathcal{H}$ & $h$ & $H=24$ \\
Representative time steps & $\mathcal{I}$ & $i$ & $I\leq T$ \\
Representative days & $\mathcal{P}$ & $p$ & $P=\nicefrac{I}{H}$ \\
Storage level time steps & $\mathcal{J}$ & $j$ & $I \leq J \leq T$ \\
\vspace{0.2cm} \\
\midrule
\midrule
Description & Sequence & Mapping & \\
\midrule
Sequence of representative time steps & $\sigma$ & $\mathcal{T}\rightarrow\mathcal{I}$ & \\
Sequence of storage level time steps & $\rho$ & $\mathcal{T}\rightarrow\mathcal{J}$ & \\
Mapping storage to representative time steps & $\vartheta$ & $\mathcal{J}\rightarrow\mathcal{I}$ & 
\end{tabular}%

  \label{tab:time_step_sets}%
\end{table}%

\subsection{Clustering representative days and representative hours}
\label{sec:clustering}
In this section, we briefly introduce key concepts of TSA with RD and RH. Given the extensive body of existing literature on various clustering algorithms \cite{Teichgraeber2022,Hoffmann2020,Nahmmacher2016,Buchholz2020}, we do not focus on the clustering algorithm. Instead, we focus on the conceptual differences between RH and RD with a focus on the sequence of the resulting representative time steps. \cref{tab:time_step_sets} summarizes the time step sets and mappings used in this paper.

The TSA receives the fully resolved time series $\tilde{x}_{a,t}$ as inputs for each attribute $a\in\mathcal{A}$ and $t\in\mathcal{T}$ with $\vert\mathcal{T}\vert=T = 8760$ hours. Typical attributes $a$ are hourly demands of electricity and heat or capacity factors of renewable technologies.
$\tilde{x}_{a,t}$ is then aggregated to $x_{a,i}$ with the representative time steps $i\in \mathcal{I}$ with $\vert\mathcal{I}\vert=I<T$. 

Important outputs of the TSA are the sequence $\sigma:\mathcal{T}\rightarrow\mathcal{I}$ of the representative time steps and the weight $w_i$ of each representative time step. The sequence $\sigma$ maps each original hour $t$ to one representative time step $i$ (\cref{tab:time_step_sets} and numbers underneath the left plot in \cref{fig:schematic_TSA_RH}).
The weight corresponds to the occurrence of the representative time step in the sequence. For example, the representative time step $i=2$ occurs 3 times in the RH aggregation in \cref{fig:schematic_TSA_RH} (in hours $t=2,t=3$, and $t=7$), hence $i=2$ has a weight of $w_2=3$.

In the case of RD, each representative time step $i$ corresponds to one hour $h \in \mathcal{H}$ with $\vert\mathcal{H}\vert=H=24$ in a representative day $p\in\mathcal{P}$ with $\vert\mathcal{P}\vert=P$. It follows: $P\times H=I$. The additional condition for RD is that the hours $h$ in each representative day $d$ are chronological. For RH, the representative time steps $i$ are independent from one another. 
Since RH aggregation does not enforce the chronology of hours, RH approximate the input time series always at least as well as RD for the same number of representative time steps $\mathcal{I}$.

All hourly-resolved variables and constraints of the optimization problem, besides the storage level variables and constraints, are aggregated by indexing them by $i\in\mathcal{I}$ instead of $t\in\mathcal{T}$. All constraints that sum up hourly values over the entire year, such as total costs or emissions, now weigh the representative time step by the corresponding weight $w_i$. The storage level constraint is the only time-coupling constraint in the presented model that requires the preservation of the time steps' chronology. Other models include other time-coupling constraints, such as ramping or minimum-up-time constraints \cite{Manco2024,Wogrin2022}, which we neglect in this work.

\subsection{Fully resolved storage constraints}
\label{sec:fully_resolved_storage_constraints}
A storage technology $s\in\mathcal{S}$ at each node $n\in\mathcal{N}$ can be charged and discharged. For the sake of readability, we omit the node index here. The difference between the storage levels $L_{s,t}$ in time step $t$ and the  previous time step $t-1$ is the net charging of the storage (charge flow $\underline{H}_{s,t}$ and discharge flow $\overline{H}_{s,t}$):
\begin{equation}
\label{eq:storage_level}
    L_{s,t}=L_{s,t-1}(1-\varphi_s)+\left(\underline{\eta}_s \underline{H}_{s,t} - \frac{\overline{H}_{s,t}}{\overline{\eta}_s}\right), \quad \forall s\in\mathcal{S},t\in\mathcal{T}\setminus\{t=1\},
\end{equation}
with the self-discharge rate $\varphi_s$, the charging efficiency $\underline{\eta}_s$, and the discharging efficiency $\overline{\eta}_s$. 
For the sake of brevity, we substitute the net charging term by $\Delta H_{s,t}$:
\begin{equation}
   \Delta H_{s,t} :=  \left(\underline{\eta}_s \underline{H}_{s,t} - \frac{\overline{H}_{s,t}}{\overline{\eta}_s}\right)
\end{equation}
To enforce storage periodicity, the storage level at $t=1$ is coupled with the last time step of the year $t=T$:
\begin{equation}
    L_{s,t=1}=L_{s,t=T}(1-\varphi_s)+\Delta H_{s,t=1}, \quad \forall s\in\mathcal{S}.
\end{equation}
$L_{s,t}$ is constrained by the energy-rated storage capacity $E_{s}$:
\begin{equation}
\label{eq:fully_resolved_storage_level_limit}
    0\leq L_{s,t}\leq E_{s}, \quad \forall s\in\mathcal{S},t\in\mathcal{T}.
\end{equation}
Furthermore, $\underline{H}_{s,t}$ and $\overline{H}_{s,t}$ are constrained by the power-rated storage capacity $C_{s}$:
\begin{equation}
    0\leq \underline{H}_{s,t} + \overline{H}_{s,t}\leq C_{s}, \quad \forall s\in\mathcal{S},t\in\mathcal{T}.\label{eq:fully_resolved_charge_limit}
\end{equation}
\cref{eq:fully_resolved_charge_limit} is a tighter formulation of the charging and discharging constraint than limiting each term individually \cite{Pozo2022}. Thus, it can help to avoid simultaneous charging and discharging. 
As with all operational variables and constraints, \cref{eq:fully_resolved_charge_limit} is reformulated in the aggregated case by replacing $t$ with $i$:
\begin{equation}
\label{eq:agg_charge_discharge_limit}
    0\leq \underline{H}_{s,i}+\overline{H}_{s,i}\leq C_{s}, \quad \forall s\in\mathcal{S},i\in\mathcal{I}.
\end{equation}
Note that in transition pathway models, such as ZEN-garden \cite{Mannhardt2025}, $E_s$ and $C_s$ are yearly resolved. However, for the sake of simplicity, we omit the year indices here.

\subsection{Proposed reduced RH storage representation}
\label{sec:met_ZEN_garden}

\begin{figure}[htbp]
    \centering
    \includegraphics[width=0.9\linewidth]{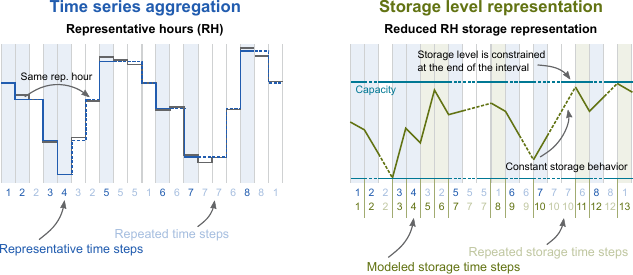}
    \caption{Schematic depiction of the time series aggregation (blue) and storage level representation (green) approaches with representative hours (RH). The fictional fully resolved time series (grey) is aggregated to 8 representative hours (numbers underneath the left plot). The storage level is represented with the presented reduced RH storage representation method for the storage level time steps (numbers underneath the right plot).} 
    \label{fig:schematic_TSA_RH}
\end{figure}

The fully resolved storage level method by \citet{Gabrielli2018} in \cref{sec:full_storage_representation_RH} - here called \textit{Full storage resolution} - has the main drawback that it requires $T$ storage level variables and $3T$ storage level constraints for each node and storage technology. As a result, the method has been found to increase computational time \cite{Blanke2022} in comparison to the method by \citet{Kotzur2018}, sometimes leading to similar solving times as the original fully resolved problem \cite{Gonzato2021}.
The novel approach presented in this study, referred to as the \textit{Proposed method} (\cref{fig:schematic_TSA_RH}), achieves a significant reduction in storage level variables and constraints compared to \textit{Full storage resolution} without compromising solution accuracy. 

The key observation to improve the \textit{Full storage resolution} method is that charging and discharging do not change over adjacent repeated time steps (\cref{fig:schematic_TSA_RH}).
As an example, the second representative hour $i=2$ in the left plot of \cref{fig:schematic_TSA_RH} is repeated twice (in hours $t=2$, and $t=3$). Hence, the net charging flows are the same for both time steps: $\Delta H_{s,i=2} =\Delta H_{s,t=2} = \Delta H_{s,t=3}$. As a result, the storage level in $t=2$ and $t=3$ can be expressed by a single variable that spans two hours. However, the same representative hour $i=2$ in hour $t=7$ cannot be represented by the same storage level time step, because the storage operation has changed in the hours $t=4$ to $t=6$. Therefore, the following rule for the storage level sequence $\rho:\mathcal{T}\rightarrow\mathcal{J}$ can be established:
\begin{equation}
    \rho(t) = \begin{cases}
    \rho(t-1), &\quad \text{if } \sigma(t)=\sigma(t-1) \\
    \rho(t-1) + 1, &\quad \text{if } \sigma(t)\neq \sigma(t-1)
    \end{cases}, \qquad \forall t\in\mathcal{T}\setminus\{t=1\}.
\end{equation}
Graphically speaking, every time the representative time step changes, a new storage level time step is added. The number of hours represented by the same storage level time step $j$ is denoted as the duration $d_j$. For example, the storage time step $j=2$ spans two hours; therefore, $d_2=2$. Since a new storage level time step $j$ is added for each change in the representative time step $i$, the mapping $\vartheta:\mathcal{J}\rightarrow\mathcal{I}$ is unique.

To derive the storage level formulation for the reduced RH storage representation, let us write the storage level equations from \cref{eq:storage_level} for the three time steps in $t=8$, $t=9$, and $t=10$:
\begin{align}
    L_{s,t=8} & = L_{s,t=7}(1-\varphi_s)+\Delta H_{s,i=5} \\
    L_{s,t=9} & = L_{s,t=8}(1-\varphi_s)+\Delta H_{s,i=5} \\
    L_{s,t=10} & = L_{s,t=9}(1-\varphi_s)+\Delta H_{s,i=5} \\
    & \nonumber\\
    \Rightarrow L_{s,t=10} & = \Bigg(\Big(L_{s,t=7}(1-\varphi_s)+\Delta H_{s,i=5}\Big)(1-\varphi_s)+\Delta H_{s,i=5}\Bigg)(1-\varphi_s)+\Delta H_{s,i=5}\\
    & = L_{s,t=7}(1-\varphi_s)^3 +\Delta H_{s,i=5}\left((1-\varphi_s)^2+(1-\varphi_s)+1\right)\\
    & = L_{s,t=7}(1-\varphi_s)^3 +\Delta H_{s,i=5}\sum_{\tilde{j}=0}^{2}(1-\varphi_s)^{\tilde{j}}
\end{align}
For all storage level time steps $j\in\mathcal{J}$, the storage level constraint is formulated as:
\begin{equation}
    \label{eq:reduced_storage_level_RH}
    L_{s,j} = L_{s,j-1}(1-\varphi_s)^{d_j}+\Delta H_{s,\vartheta(j)}\sum_{\tilde{j}=0}^{d_j-1}(1-\varphi_s)^{\tilde{j}}, \quad \forall s\in\mathcal{S},j\in\mathcal{J}\setminus\{j=1\}.
\end{equation}
The other constraints are formulated equivalently:
\begin{align}
   & L_{s,j=1}=L_{s,j=J}(1-\varphi_s)^{d_{j=1}}+\Delta H_{s,\vartheta(j=1)}\sum_{\tilde{j}=0}^{d_{j=1}-1}(1-\varphi_s)^{\tilde{j}}, \quad && \forall s\in\mathcal{S}, &\\
   & 0\leq L_{s,j} \leq E_s, \quad && \forall s\in\mathcal{S},j\in\mathcal{J}.& \label{eq:reduced_storage_level_limit_RH}
\end{align}
We can show that $L_{s,j}$ is monotonic over all hours $t$ in $j$. Monotony implies that \cref {eq:reduced_storage_level_limit_RH} provides an upper and lower limit to the fully resolved storage level $L_{s,t}$ for all $t\in\mathcal{T}$. Consider \cref{eq:reduced_storage_level_RH} for one storage level time step $j\in\mathcal{J}$, during which $\Delta H_{s,\vartheta(j)}$ is constant. The storage level $L_{\hat{t}}$ for one generic storage time step $j$ and storage technology $s$ for the intermediate time steps $\hat{t}\in[1,d_j]$ follows as:
\begin{equation}
\label{eq:storage_level_simpl_mon}
    L_{\hat{t}} = L_0(1-\varphi)^{\hat{t}}+\Delta H\sum_{\tilde{t}=0}^{\hat{t}-1}(1-\varphi)^{\tilde{t}},
\end{equation}
with the storage level at the end of the previous storage level time step $L_0 = L_{j-1}$. Without self-discharge ($\varphi=0$), it follows for the derivative of $ L_{\hat{t}}$:
\begin{equation}
    \dv{L_{\hat{t}}}{\hat{t}}=\Delta H.
\end{equation}
Therefore, $L_{\hat{t}}$ is monotonic over $\hat{t}$ for $\varphi=0$. 

For $0<\varphi<1$, $\sum_{\tilde{t}=0}^{\hat{t}-1}(1-\varphi)^{\tilde{t}}$ is reformulated as the partial geometric series:
\begin{equation}
    L_{\hat{t}} = L_0(1-\varphi)^{\hat{t}}+\Delta H\frac{1-(1-\varphi)^{\hat{t}}}{\varphi}.
\end{equation}
It follows for $\dv*{L_{\hat{t}}}{\hat{t}}$:
\begin{equation}
    \dv{L_{\hat{t}}}{\hat{t}} = \underbrace{\left(L_0-\frac{\Delta H}{\varphi}\right) \ln(1-\varphi)}_{=\text{ constant } \forall \hat{t}\in[1,d_j]}(1-\varphi)^{\hat{t}}.
\end{equation}
With $(1-\varphi)^{\hat{t}}>0$, it follows that \cref{eq:storage_level_simpl_mon} is monotonic for $0<\varphi<1$.
Thus, the storage level constraints \cref{eq:reduced_storage_level_limit_RH} apply to the entire time interval $\mathcal{T}$, even though not every hour $t\in\mathcal{T}$ is explicitly modeled. The storage level evolution for each $t\in\mathcal{T}$ can be reconstructed ex-post.

\section{Case study of a net-zero European electricity and heating system}
\label{sec:case_study}
The presented storage representation method from \cref{sec:methods} is applied to a case study of the European electricity and heating system to test the method's accuracy and performance in spatially resolved large-scale energy systems. We minimize total annual system costs in the objective function.

Based on the case study from \citet{Mannhardt2024a}, we conduct a greenfield optimization of a carbon-neutral energy system in 2050, in which long-term storage is relevant to constantly provide emission-free electricity and heat. 
The case study contains 28 nodes, corresponding to European countries, and encompasses 31 conversion technologies, 3 transport technologies, and 4 storage technologies (batteries, hydrogen storage, pumped hydro storage, and natural gas storage). Besides 15 electricity generation technologies and 13 heat generation technologies, the case study considers carbon storage, LNG terminals, and industrial gas consumers.

The time-dependent input data are the electricity and heating demand, the capacity factor of solar photovoltaics (PV), onshore wind, offshore wind, reservoir hydropower, and run-of-river hydropower, and the conversion factor of residential and district heating heat pumps. Furthermore, each heating technology has a time-dependent capacity factor, which is calculated as the share of the hourly heat demand of the peak heat demand. 
The motivation and derivation of the capacity factors for heating technologies are described in detail in the Supplementary Information of \cite{Mannhardt2023}. 
 
\section{Results}
\label{sec:results}
We benchmark the proposed reduced RH storage representation, \textit{Proposed method}, with the four established storage representation methods: The first two methods, i.e., \textit{Superposition} and \textit{MinMax} \cite{Kotzur2018}, utilize RD and differ in whether the storage level limit is constrained in each hour or once per day, respectively. See \cref{sec:storage_representation_RD} for a detailed mathematical description of these methods. The third method, i.e., \textit{Full storage resolution}, utilizes RH and formulates the storage level for each hour of the year \cite{Gabrielli2018}. See \cref{sec:full_storage_representation_RH} for a detailed description of \textit{Full storage resolution}. The fourth method, i.e., \textit{Chrono}, aggregates the input time series to CRH and represents the storage time index with the same chronological time index \cite{Pineda2018}. See \cref{sec:storage_representation_CRH} for a detailed description of \textit{Chrono}.
Note that \textit{Chrono} is a special case of \textit{Proposed method} in \cref{sec:met_ZEN_garden}, and the same implementation can be used.

First, we compare the error in the objective value, i.e., total system cost, with the solving time. Afterwards, we investigate the error in conversion and storage capacity, and annual storage cycles. The benchmarking is completed by assessing the robustness of the proposed methods in accurately optimizing variations of the investigated case study.

\begin{table}[htbp]
    \small
  \centering
  \caption{Summary of benchmarked storage representation methods, the utilized time series aggregation method (TSA), and how the storage level is represented and constrained.}
\begin{tabular}{lllll}
Method & TSA   & Storage representation & Storage level limit & Reference \\
\midrule
\textit{Proposed method} & RH    & Reduced & Every storage time step & This paper \\
\textit{Superposition} & RD    & Superposition & Every hour & \citet{Kotzur2018} \\
\textit{MinMax} & RD    & Superposition & Twice a day & \citet{Kotzur2018} \\
\textit{Full storage resolution} & RH    & Full resolution & Every hour & \citet{Gabrielli2018} \\
\textit{Chrono} & CRH   & Reduced & Every storage time step & \citet{Pineda2018} \\
\end{tabular}%

  \label{tab:benchmark_methods}%
\end{table}%

All storage representation methods are implemented in the open-source ESOM ZEN-garden \cite{Mannhardt2025} and are available in the \href{https://github.com/jacob-mannhardt/ZEN-garden-fork/tree/publication_storage_representation}{Github repository}. The documentation of the ESOM is available online and includes the mathematical formulation of ZEN-garden \cite{Mannhardt2025a}. The input data and results for all models and methods are found in \cite{Mannhardt2025c}.

We optimize the case study presented in \cref{sec:case_study} for each storage representation method from \cref{sec:methods} and 10 different time series aggregation settings. The strongest aggregation is 24 representative hours (for RH and CRH) or 1 representative day (for RD). The subsequent aggregations are 48 RH, 96 RH, 192 RH, 384 RH, 768 RH, 1536 RH, 3072 RH, 6144 RH, and finally full resolution with 8760 RH. 

In this work, we use the time series aggregation module \textit{tsam} 2.1.0 \cite{Hoffmann2022}.
Since the focus of this work is on the storage level representation and not on the clustering algorithms for time series aggregation, we run all scenarios with hierarchical clustering using mean representation, except for CRH by \citet{Pineda2018}, which is clustered by adjacent periods. \citet{Hoffmann2021} show that the clustering algorithm for non-chronological time steps has little impact on the comparison between RD and RH.
\subsection{Objective value versus solving time}
\begin{figure}[tbph]
    \centering
    \includegraphics[width=0.9\linewidth]{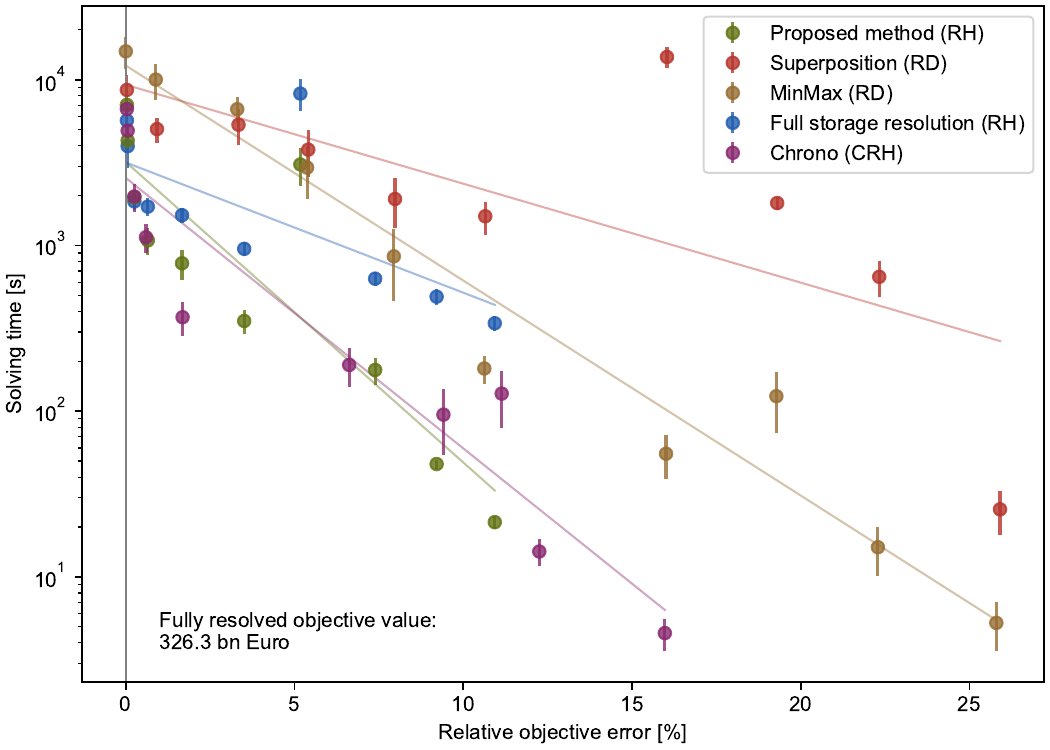}
    \caption{Solving time versus relative objective error for the five investigated storage representation methods. Every method is repeated five times for each time step configuration. The error bar indicates the 95\% confidence interval in solving time. Note that the y-axis is shown logarithmically. The linear lines show an exponential regression of the mean solving time and objective value.} 
    \label{fig:time_vs_objective}
\end{figure}

The goal of TSA and storage representation is to approximate the fully resolved solution as closely as possible while reducing the solving time as much as possible. \cref{fig:time_vs_objective} compares the objective value relative to the fully resolved solution (326.3 bn Euro) with the solving time for the five investigated storage representation methods. 
Since all approaches precisely track the storage level, the objective value is only determined by the TSA. Hence, the two RD approaches, \textit{Superposition} and \textit{MinMax}, and the two RH approaches, \textit{Full storage resolution} and \textit{Proposed method}, are vertically aligned in \cref{fig:time_vs_objective}, respectively.

\textit{MinMax} shows shorter solving times than  \textit{Superposition} for strong aggregation, but leads to longer solving times close to the fully resolved solution. \textit{Full storage resolution} has the highest solution times for the most strongly aggregated model, but the increase for higher resolutions is the lowest of all representation methods. The variable and constraint reduction of \textit{Proposed method} strongly reduces solving time in comparison to \textit{Full storage resolution}. For example, the most aggregated model (24 RH) solves on average 20 times faster for \textit{Proposed method} than for \textit{Full storage resolution}. \textit{Chrono} shows the shortest solving times at the strongest aggregation across all models.

For the strongest aggregation (24 RH or 1 RD), the RD approaches show a relative objective error of 25.9\%, whereas the RH approaches achieve a relative objective error of 10.9\%. The CRH approach shows a relative objective error of 16\%. The reason for the discrepancy in the objective value between RD and RH reflects the higher degree of freedom for RH in the TSA. With one representative day, the RD approaches cannot capture any seasonal behavior, whereas the RH approaches can spread the 24 time steps throughout the year and thereby capture seasonality. At 16 RD (384 RH), the RD approaches show a similar objective value as the RH approaches at 24 RH. In detail, the RH approaches show a smaller objective error than the RD approaches for all time step configurations. The CRH approach shows slightly increased objective errors compared to the RH approaches throughout.

The differences in solving time can be attributed to the size of the optimization problem, specifically the number of storage level variables and constraints (\cref{fig:num_cons_vars}). In general, the more storage level variables and constraints are formulated, the larger is the resulting optimization problem, and, therefore, the longer it usually takes to solve \cite{Spielman2003,Wright1997}. 
\textit{Full storage resolution} has a constant number of storage level variables ($T=8760$ per technology and node) and constraints ($3T=26 280$ per technology and node). Therefore, \textit{Full storage resolution} shows the highest solving time at the strongest aggregation (\cref{fig:time_vs_objective}). The RD methods \textit{Superposition} and \textit{MinMax} significantly reduce the number of storage level variables compared to \textit{Full storage resolution}, which reduces the solving time for strong aggregation. \textit{MinMax} further reduces the number of storage level constraints, leading to mostly lower solving times compared to \textit{Superposition}. \textit{Proposed method} and \textit{Chrono} show significantly fewer storage level constraints and variables than \textit{Full storage resolution}. 

\begin{figure}[tbph]
    \centering
    \includegraphics[width=0.75\linewidth]{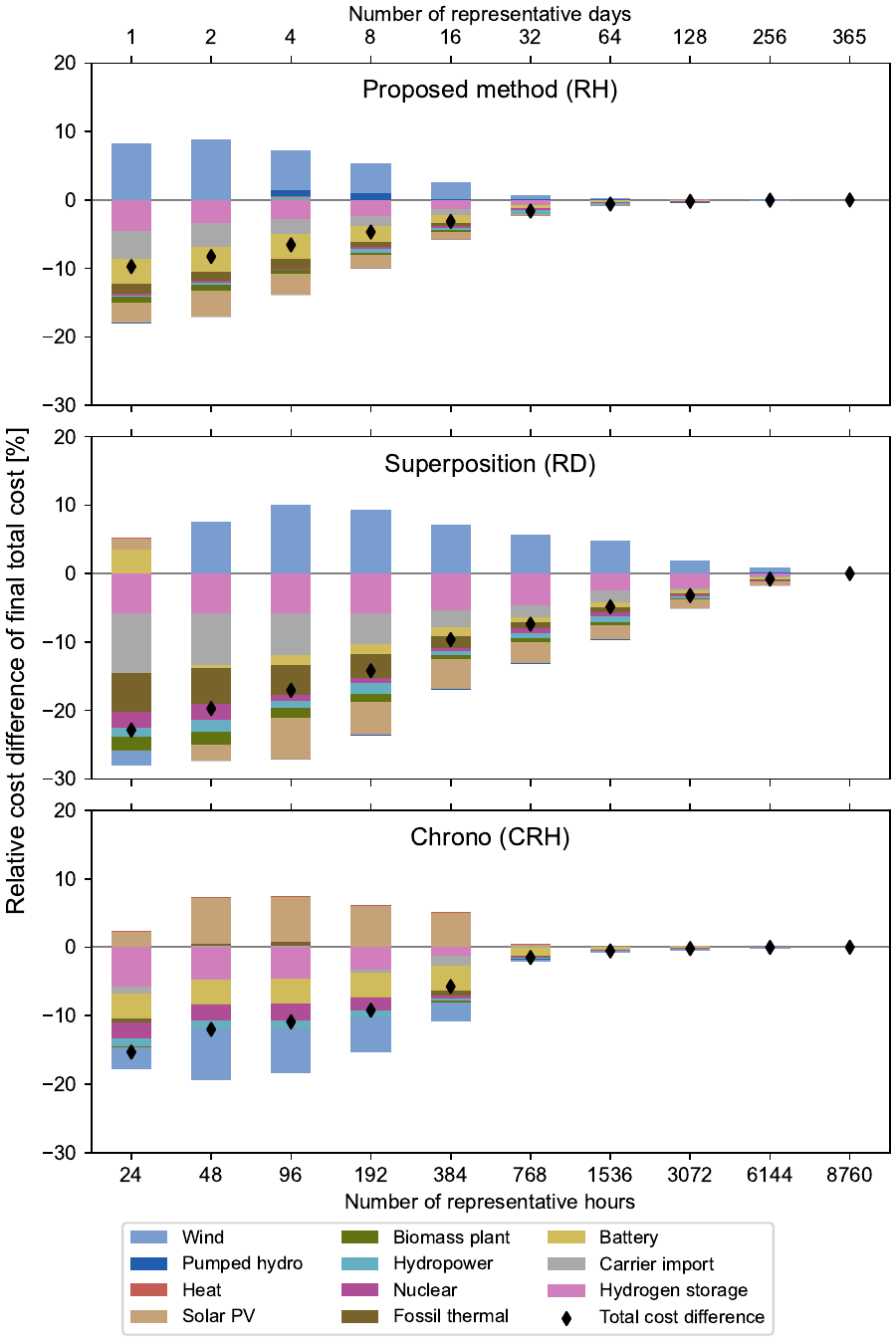}
    \caption{Disaggregated cost error for \textit{Proposed method} (RH), \textit{Superposition} (RD), and \textit{Chrono} (CRH), relative to fully resolved total cost (326.3 bn Euro). \textit{Full storage resolution} and \textit{MinMax} are not shown, as they have the same system design as \textit{Proposed method} and \textit{Superposition}, respectively.} 
    \label{fig:cost_error}
\end{figure}

However, \cref{fig:cost_error} shows that a similar objective value does not imply that the system design is similar across solutions. \textit{Proposed method} and \textit{Chrono} show a similar total cost difference, but \textit{Proposed method} relies more heavily on wind and shows slightly smaller technology-specific errors than \textit{Chrono}, except for in the first two time step configurations. For the strongest aggregation of 24 RH (1 RD), \textit{Chrono} shows the lowest technology-specific cost errors. Then, for higher resolutions, the model shows a strongly positive error for solar PV and a strongly negative error for wind, which persists up until 768 RH. \textit{Proposed method} and \textit{Chrono} continuously reduce the cost error and show virtually no positive cost error from 768 RH onward. For all time step configurations, \textit{Superposition} shows reduced carrier imports and less investment in storage technologies, especially hydrogen storage. Lower investments in solar PV are offset by increased onshore wind, except for 1 RD. The persistent error in \textit{Superposition} leads to the highest errors across all investigated approaches.

\subsection{Approximation accuracy of conversion and storage capacities}
\begin{figure}[htbp]
    \centering
    \includegraphics[width=0.75\linewidth]{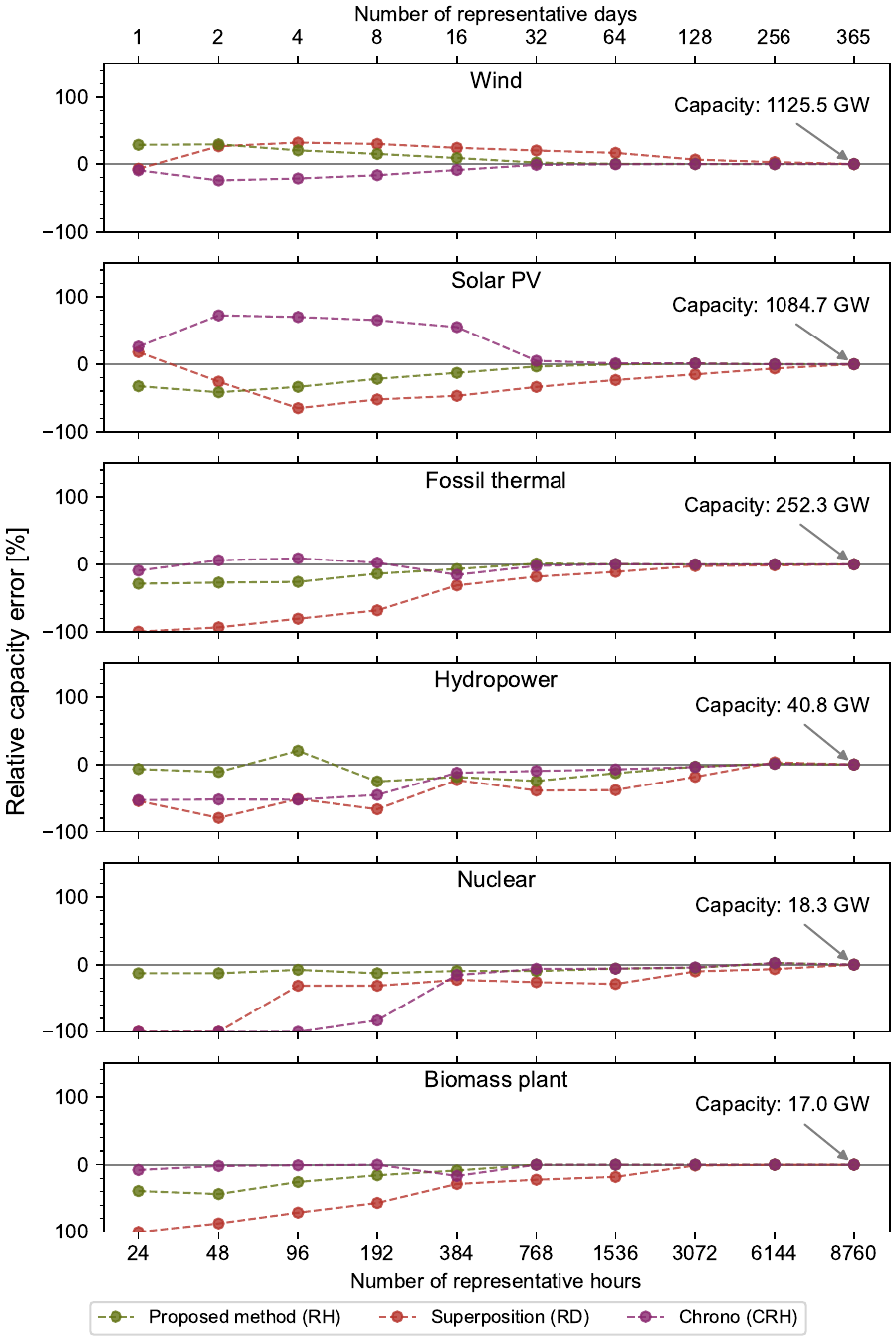}
    \caption{Relative error in electricity generation capacity for \textit{Superposition} (RD), \textit{Proposed method} (RH), and \textit{Chrono} (CRH), relative to fully resolved capacities. The technologies are ranked by total capacity. \textit{Full storage resolution} and \textit{MinMax} are not shown, as they have the same system design as \textit{Proposed method} and \textit{Superposition}, respectively.} 
    \label{fig:capacity_error_conversion}
\end{figure}

\cref{fig:capacity_error_conversion} shows the relative error in electricity generation capacity for \textit{Proposed method}, \textit{Superposition}, and \textit{Chrono}. For most approaches, a stronger aggregation results in less capacity of dispatchable technologies, such as fossil thermal, biomass, and nuclear power plants. Since the clusters are represented by their mean value, the time series are smoothed, which reduces the need for dispatchable technologies. Wind and solar PV are the most prevalent technologies with 1126 GW and 1085 GW under full resolution, respectively. Fossil thermal power plants, especially natural gas turbines, contribute 252 GW. All other technologies together account for less than 100 GW.

Overall, \textit{Proposed method} shows the closest approximation of electricity generation capacities over all time step configurations, especially for the strongest aggregation. The most dominant technologies, wind and solar PV, are approximated without significant error at 768 RH, whereas \textit{Superposition} requires a full resolution to achieve the same accuracy. In general, \textit{Proposed method} overestimates wind capacity and underestimates solar PV capacity; however, the relative error always remains below 40\% (maximum error for solar PV for 48 RH).

Even though \textit{Proposed method} and \textit{Chrono} both cluster RH and use the same storage level representation method, the capacity installations behave inversely: \textit{Chrono} overestimates the solar PV capacity (up to +73\% relative error), and in turn shows lower wind capacity (up to -21\% relative error). The reason for the overreliance of \textit{Chrono} on solar PV is that the chronological clustering leads to long periods with the same representative time step. As a result, the capacity factor of solar PV is smoothed to the same value over many days, and the daily intermittence is removed. From 1536 RH onward, all capacities are approximated closely.

\textit{Superposition} shows the highest error in dispatchable capacity, and especially the capacity of fossil thermal, biomass, and nuclear power plants is strongly underestimated at strong aggregations. The step from 1 RD (24 RH) to 2 RD (48 RH) leads to a jump in the capacity approximation of wind and solar PV: For 1 RD, only daily behavior can be modeled, and, therefore, the optimizer relies more on solar PV (overestimation of +18\%). For 2 RD, some differences between days are modeled, and thus, one of the two RD has low solar irradiance, which reduces the value of solar PV. Hence, \textit{Superposition} shows a negative capacity error for solar PV (up to -65\%) from 2 RD onward and installs more wind capacity (up to +32\%).

\begin{figure}[htbp]
    \centering
    \includegraphics[width=0.75\linewidth]{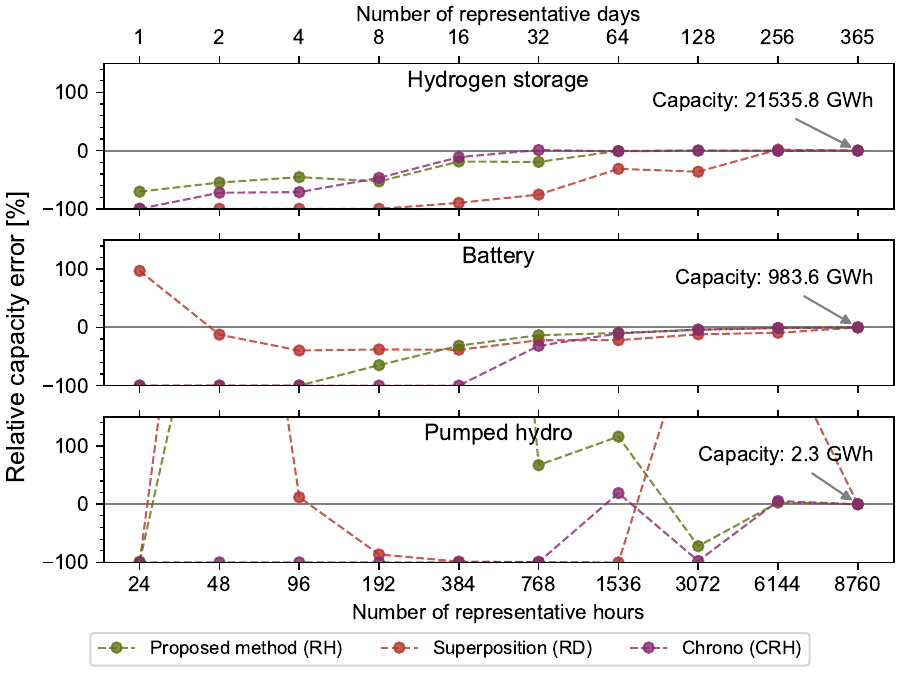}
    \caption{Relative error in electricity storage capacity for \textit{Superposition} (RD), \textit{Proposed method} (RH), and \textit{Chrono} (CRH), relative to fully resolved capacities. The technologies are ranked by total capacity. \textit{Full storage resolution} and \textit{MinMax} are not shown, as they have the same system design as \textit{Proposed method} and \textit{Superposition}, respectively.} 
    \label{fig:capacity_error_storage}
\end{figure}

The different TSA methods lead to divergent results for short- and long-term storage investment (\cref{fig:capacity_error_storage}). Hydrogen storage is a promising long-term storage technology thanks to low energy-rated investment costs, high power-rated investment costs, low self-discharge, and low round-trip efficiency. Batteries, on the other hand, are best suited for short-term storage, since they have comparably low power-rated investment costs, high energy-rated investment costs, higher self-discharge, and high round-trip efficiency. Pumped hydro storage is generally considered a medium-term storage option between batteries and hydrogen storage. \cref{fig:capacity_error_storage} shows that with full resolution 21.5 TWh of hydrogen storage and 984 GWh of battery storage are built. Pumped hydro storage plays a minor role with 2.3 GWh. Note that we do not consider existing capacities in this study; therefore, the existing pumped hydro storage is neglected.

\textit{Proposed method} shows a passable approximation of hydrogen storage at a low resolution. At 24 RH, hydrogen storage capacities are around 70\% below the fully resolved solution, whereas \textit{Superposition} and \textit{Chrono} show no hydrogen storage capacities at all. For higher resolutions, \textit{Proposed method} and \textit{Chrono} continuously decrease the capacity error, while the results with \textit{Superposition} do not show any significant hydrogen storage capacities until 384 RH (16 RD). 

On the other hand, \textit{Superposition} overestimates battery capacity at 1 RD by 100\% due to the purely daily behavior, which imposes a positive bias on battery storage. \textit{Proposed method} and \textit{Chrono} show no battery capacity for 24, 48, and 96 RH, after which battery capacity is increased under \textit{Proposed method} and shows the same capacity as in \textit{Superposition} for 384 RH (-31\% relative error). Hence, \textit{Proposed method} shows a good approximation of both long-term and short-term storage from around 384 RH. 
\textit{Chrono} does not show any significant battery capacity until 768 RH, and therefore underestimates the overall need for short-term storage. In general, the time series smoothing in  \textit{Chrono} leads to a reduced reliance on both short- and long-term storage.
Pumped hydro storage is in the twilight between playing a significant role and being outcompeted by batteries and hydrogen storage, so its relative error fluctuates strongly across all approaches. 

\begin{figure}[htbp]
    \centering
    \includegraphics[width=0.75\linewidth]{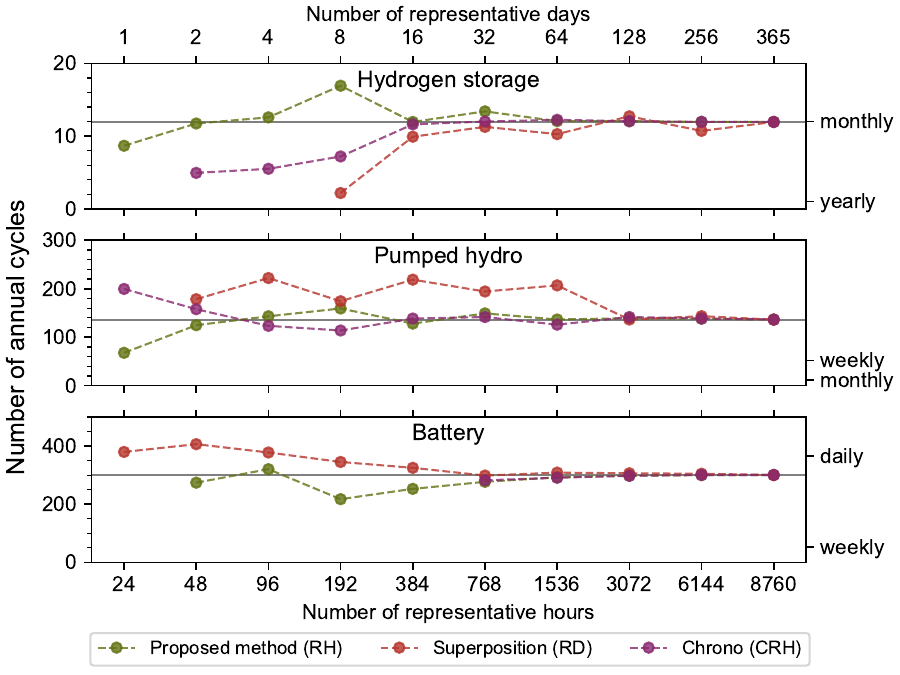}
    \caption{Number of annual cycles for\textit{Proposed method} (RH),  \textit{Superposition} (RD), and \textit{Chrono} (CRH). The technologies are ranked by an increasing number of cycles. The right y-axis indicates the corresponding daily, weekly, monthly, or yearly pattern.
    \textit{Full storage resolution} and \textit{MinMax} are not shown, as they have the same system design as \textit{Proposed method} and \textit{Superposition}, respectively. No value is shown for capacity values below 0.01\% of the fully resolved capacity (\cref{fig:capacity_error_storage}).} 
    \label{fig:storage_cycles}
\end{figure}

The TSA method impacts the role that the storage technologies play in the energy system (\cref{fig:storage_cycles}). In the fully resolved energy system, hydrogen storage is a long-term storage technology and shows monthly cycles (12 cycles per year). Pumped hydro and battery storage operate more as short-term storage with 136 (half-weekly) and 300 (almost daily) cycles per year, respectively. \textit{Proposed method} approximates the number of cycles well for the entire range of time step configurations. \textit{Superposition} relies significantly less on hydrogen storage (\cref{fig:capacity_error_storage}) but shows an accurate number of cycles once installed. In return, \textit{Superposition} shows daily battery cycles for strong aggregation (365 cycles) and thereby generally overestimates battery cycles.
The long time step durations in \textit{Chrono} lead to fewer storage cycles for hydrogen storage, as the technology is only used for multi-monthly storage. Once batteries are installed at 768 RH, they show an accurate number of cycles, which highlights that batteries are only installed if they can operate on a quasi-daily cycle.

The solution accuracy is further highlighted by the price duration curve of the electricity shadow prices (\cref{fig:snapshot_TSA_electricity_price_LDS}). \textit{Proposed method} and \textit{Chrono} approximate the electricity price well for most TSA settings, whereas \textit{Superposition} shows significantly lower maximum and higher minimum electricity prices. \textit{Chrono} shows a similar behavior to \textit{Proposed method} at the highest prices, but the price drops to almost zero for half of the year for 24 RH.

\subsection{Robustness of methods for various large-scale energy system models}

\begin{figure}[htbp]
    \centering
    \includegraphics[width=0.75\linewidth]{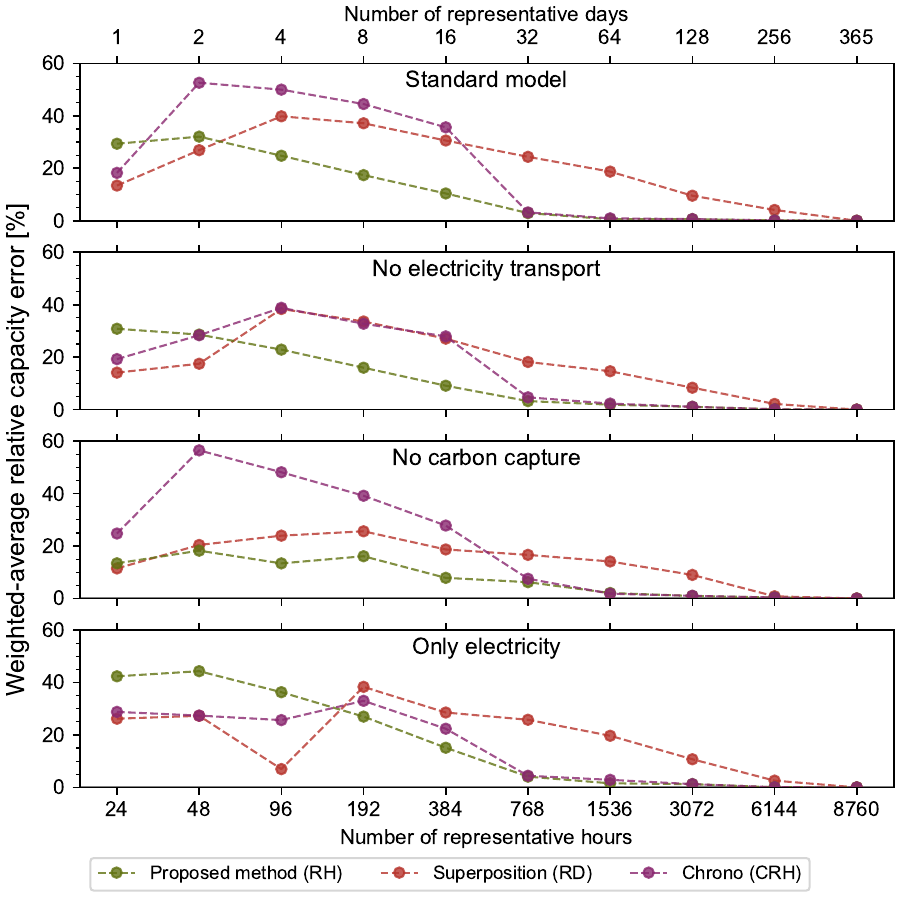}
    \caption{Weighted-average relative absolute capacity error of all electricity generation technologies in four large-scale energy system models for \textit{Superposition} (RD), \textit{Proposed method} (RH), and \textit{Chrono} (CRH). \textit{Full storage resolution} and \textit{MinMax} are not shown, as they have the same system design as \textit{Proposed method} and \textit{Superposition}, respectively.} 
    \label{fig:average_capa_error}
\end{figure}

\cref{fig:average_capa_error} extends the previous results by investigating three additional large-scale energy system models, all based on the case study presented in \cref{sec:case_study}. The results show the robustness of the storage representation methods and highlight cases in which the approaches perform particularly well or badly.
First, electricity transport via power lines across nodes is disabled; therefore, each country must produce its own electricity. Second, carbon capture is removed, which requires the energy system to provide carbon-neutral electricity and heat at every moment without carbon capture. Third, the heating sector is removed to investigate an energy system with predominantly daily time series.

\cref{fig:average_capa_error} shows the weighted-average relative absolute capacity error of all electricity generation technologies for the four case studies. First, the relative capacity error is calculated for each technology, analogously to \cref{fig:capacity_error_conversion}. Then, the absolute values of the relative error are weighted by the capacity and averaged. Across all four models, \textit{Chrono} shows a high error for strong aggregation, except for 24 RH, but an accurate representation for many representative hours. \textit{Superposition}, on the other hand, shows higher accuracy than \textit{Chrono} for a strong aggregation but retains a significant error at higher resolutions.

For the standard model and the models without either electricity transport or carbon capture, \textit{Proposed method} shows a lower error than both \textit{Superposition} and \textit{Chrono} for medium-low resolutions (96 RH to 768 RH). At higher resolutions, the error is indistinguishable from \textit{Chrono}. When only considering the electricity sector, the daily pattern of the time series supports the use of RD, and therefore, \textit{Superposition} performs best for low resolutions.

Overall, \cref{fig:average_capa_error} shows a robust ability of \textit{Proposed method} to provide a satisfiable accuracy at resolutions of around 100 to 500 RH across models. An average capacity error of below 10\% is achieved at 384 RH, except for the ``only electricity'' model. \textit{Chrono} shows a higher capacity error at low resolutions, because of the inability to capture daily behavior. \textit{Superposition} needs significantly more representative time steps to reduce the capacity error as much as \textit{Proposed method} and \textit{Chrono}.
\section{Discussion}
\label{sec:disc}
Large-scale energy system models with many regions, technologies, interconnected sectors, and a transition pathway require complexity reduction methods to be computationally feasible. Most commonly, time series aggregation reduces the intra-annual time complexity. However, the time coupling constraints of storage technologies complicate the reduction and often lead to strong misrepresentations of storage and conversion capacities.

This paper presents a novel storage representation method to combine the accuracy of clustering representative hours (RH) with the computational efficiency of existing storage representation methods using representative days (RD). Clustering RH can approximate the input data significantly more accurately than RD for the same number of representative time steps. However, until now, the efficient representation of short- and long-term storage with self-discharge at strong time series aggregation has not been achieved using RH. 


We benchmark the presented storage representation approach 
with two approaches relying on RD \cite{Kotzur2018}, one RH method with a fully resolved storage level \cite{Gabrielli2018}, and one method with chronological RH clustering (CRH) \cite{Pineda2018}.

Our results show that the RH and CRH methods approximate the objective value, disaggregated cost terms, and capacities more accurately than the RD approaches for the same number of representative time steps. The CRH method is a special case of the roposed method, and the results are similar between the two representation methods.

In detail, the objective error of the proposed method is approximately 10\% for the strongest aggregation to 24 RH. The capacities of the RH approaches converge mostly smoothly to the fully resolved value, whereas the CRH method shows slightly higher errors for a low number of representative hours, and the RD methods show persistent errors for higher resolutions.
The RH method with a fully resolved storage level shows the highest solving time, which is reduced by up to 20-fold in the proposed method, while retaining the same solution quality. 
Overall, the proposed method and the CRH method show the most accurate approximation at the lowest computational burden. Since CRH clusters chronological time steps, it underestimates daily variability, leading to a significant misrepresentation of solar PV and battery capacities and operations. 
Furthermore, the CRH method underestimates the number of cycles in long-term storage technologies.
The proposed method strikes a good balance between conserving daily and long-term patterns.

Our multi-case-study investigation reveals a robust behavior of the benchmarked approaches across different model configurations. The most notable difference is observed when optimizing the electricity sector in isolation, which shows mainly daily time series patterns. In this case, the proposed method showed a higher electricity capacity error than the other methods for low resolutions. This observation warrants further comparisons of the methods for different energy system designs, including smaller scales, such as multi-energy systems \cite{Gabrielli2018,Srinivasan2025}, or non-European case studies \cite{Parker2017}.


In this work, we do not discuss the impact of different time series clustering methods on the applicability of the storage representation approaches. \citet{Hoffmann2021} have shown that the clustering algorithm plays a subordinate role in comparison to selecting RH or RD. However, investigating various clustering algorithms could reveal optimal pairings between clustering and storage representation methods. For example, a multi-objective clustering method that, in addition to minimizing the mean square error of the clustering, pays special attention to clusters of adjacent hours without being strictly chronological would support the developed storage representation method. 

Extreme period behavior becomes increasingly important with a higher penetration of renewable technologies \cite{Bahl2018,Teichgraeber2020,DeMarco2025}. While we exclude their discussion in this work, it is valuable to understand how the different time series aggregation and storage level representation approaches manage to account for extreme period behavior. Intuitively, RH may be better positioned than RD to include single extreme hours in the sequence of representative time steps.

\subsection{Conclusions}
The proposed method solves the trade-off between accuracy and computational tractability by:
\begin{enumerate}
    \item Clustering RH to approximate the input time series accurately with a low number of representative time steps,
    \item Reducing the number of required storage level variables and constraints by combining adjacent storage level time steps that are represented by the same RH.
\end{enumerate}
The presented approach is particularly suited for strong aggregation to around 100 to 500 representative time steps, which is required for multi-year transition pathway optimizations of spatially resolved, sector-coupled energy system models. 

The majority of energy system optimization models with long-term storage utilize RD and the superposition method by \citet{Kotzur2018}. 
We demonstrate that it is possible and in many cases preferable to use RH while retaining an efficient and accurate storage level formulation. 

\section*{Acknowledgments}
We thank Dr. Leonard Göke and Dr. Paolo Gabrielli for their valuable feedback on the manuscript.

The research published in this report was carried out with the support of the Swiss Federal Office of Energy (SFOE) as part of the SWEET PATHFNDR project. The authors bear sole responsibility for the conclusions and the results.
\newpage
\appendix

\section{Mathematical description of benchmarked storage representation methods}

\begin{table}[htbp]
  \centering
  \caption{Summary of time step sets and mappings for representative days.}
\begin{tabular}{llll}
Description & Set & Index & Dimensionality  \\
\midrule
Days per year & $\mathcal{D}$ & $d$ & $D=\nicefrac{T}{H}=365$ \\
\vspace{0.2cm} \\
\midrule
Description & Sequence & Mapping & \\
\midrule
Sequence of representative days & $\psi$ & $\mathcal{D}\rightarrow\mathcal{P}$ & \\
\end{tabular}%

  \label{tab:time_step_sets_appendix}%
\end{table}%

\subsection{Storage representation approaches for representative days}
\label{sec:storage_representation_RD}

\begin{figure}[tbph]
    \centering
    \includegraphics[width=0.9\linewidth]{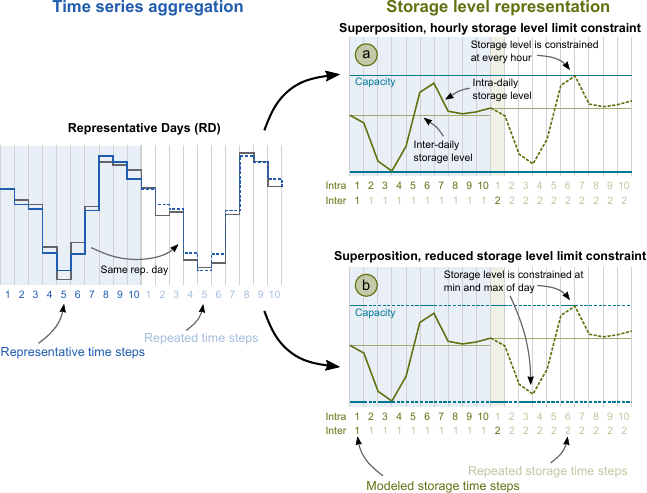}
    \caption{Schematic depiction of the time series aggregation (blue) and storage level representation (green) approaches with representative days (RD). The fictional fully resolved time series (grey) is aggregated to one representative day with 10 chronological hours per day. The storage level is represented by the superposition method of \cite{Kotzur2018} with hourly (\textit{Superposition}, panel a) or reduced storage level limit constraint (\textit{MinMax}, panel b).} 
    \label{fig:schematic_TSA_RD}
\end{figure}

The most common storage representation approach for representative days (RD) by \citet{Kotzur2018} superposes an inter-daily storage level $L^\mathrm{inter}_{s,d}$ for each day $d\in\mathcal{D}$ $(\vert\mathcal{D}\vert =D=365)$ with an intra-daily storage level $L^\mathrm{intra}_{s,p,h}$ for each hour $h$ of a representative day $p\in\mathcal{P},~\vert\mathcal{P}\vert=P\leq D$.

The constraint for $L^\mathrm{intra}_{s,p,h}$ is analogous to the unaggregated storage level equation in \cref{eq:storage_level}:
\begin{equation}
\label{eq:intra_level_kotzur}
    L^\mathrm{intra}_{s,p,h+1} = L^\mathrm{intra}_{s,p,h}(1-\varphi_s)+\Delta H_{s,p,h}, \quad \forall s\in\mathcal{S},p\in\mathcal{P},h\in\mathcal{H}\setminus\{h=H\}.
\end{equation}
Note that we shift the hour index $h$ in comparison to the formulation in \cref{eq:storage_level} to follow the original formulation in \cite{Kotzur2018}. Both formulations are correct, but calculate the storage level either after (\cref{eq:storage_level}) or before (\cref{eq:intra_level_kotzur}) the current hour $h$.
$L^\mathrm{intra}_{s,p,h}$ at the beginning of the day $h=1$ is zero:
\begin{equation}
    L^\mathrm{intra}_{s,p,h=1}=0, \quad \forall s\in\mathcal{S},p\in\mathcal{P}.
\end{equation}
$L^\mathrm{inter}_{s,d}$ between days $d$ is coupled through $L^\mathrm{intra}_{s,p,h}$ in the last hour of the representative day $p$ with the mapping $\psi:\mathcal{D}\rightarrow \mathcal{P}$ (\cref{tab:time_step_sets_appendix}):
\begin{equation}
\label{eq:inter_intra_coupling}
    L^\mathrm{inter}_{s,d+1} = L^\mathrm{inter}_{s,d}(1-\varphi_s)^{H}+L^\mathrm{intra}_{s,\psi(d),H+1},\quad \forall s\in\mathcal{S},d\in\mathcal{D}.
\end{equation}
We slightly improve the original formulation by substituting $L^\mathrm{intra}_{s,\psi(d),H+1}$ by the corresponding formulation from \cref{eq:intra_level_kotzur}, which reduces the number of required variables:
\begin{equation}
    L^\mathrm{inter}_{s,d+1} = L^\mathrm{inter}_{s,d}(1-\varphi_s)^{H}+L^\mathrm{intra}_{s,\psi(d),H}(1-\varphi_s)+\Delta H_{s,\psi(d),H},\quad \forall s\in\mathcal{S},d\in\mathcal{D}\setminus \{d=D\}.
\end{equation}
As before, the periodicity constraint is formulated by connecting the first and the last day of the year:
\begin{equation}
    L^\mathrm{inter}_{s,1} = L^\mathrm{inter}_{s,D}(1-\varphi_s)^{H}+L^\mathrm{intra}_{s,\psi(D),H}(1-\varphi_s)+\Delta H_{s,\psi(D),H},\quad \forall s\in\mathcal{S}.
\end{equation}
Note that the storage level is not explicitly implemented in the optimization problem but can be calculated ex-post.
\subsubsection{Superposition, hourly storage level limit constraint}
\label{sec:met_kot_orig}
In the original publication by \citet{Kotzur2018}, the storage level is constrained by the capacity for each hour $t\in\mathcal{T}$ with $T=D\times H$ (panel a in \cref{fig:schematic_TSA_RD}):
\begin{equation}
\label{eq:kotzur_hourly_storage_level_limit}
    0\leq L^\mathrm{inter}_{s,d}(1-\varphi_s)^h + L^\mathrm{intra}_{s,\psi(d),h} \leq E_s, \quad \forall s\in\mathcal{S},d\in\mathcal{D},h\in\mathcal{H}.
\end{equation}
Note that we fixed the error in the original formulation by changing the term $(1-\varphi_s)^H$ to $(1-\varphi_s)^h$, as observed by \cite{Blanke2022}.

\subsubsection{Superposition, reduced storage level limit constraint}
\citet{Kotzur2018} note in the appendix of the original publication that the number of constraints necessary for \cref{eq:kotzur_hourly_storage_level_limit} can be reduced by limiting the storage level by the maximum and minimum intra-daily storage level $L^\mathrm{intra,max}_{s,p}$ and $L^\mathrm{intra,min}_{s,p}$ (panel b in \cref{fig:schematic_TSA_RD}):
\begin{align}
    L^\mathrm{intra,max}_{s,p} &\geq L^\mathrm{intra}_{s,p,h}, \quad \forall s\in\mathcal{S},p\in\mathcal{P},h\in\mathcal{H}, \\
    L^\mathrm{intra,min}_{s,p} &\leq L^\mathrm{intra}_{s,p,h}, \quad \forall s\in\mathcal{S},p\in\mathcal{P},h\in\mathcal{H}.
\end{align}
Then, the storage level is constrained for every day $d\in\mathcal{D}$ as:
\begin{align}
    & L^\mathrm{inter}_{s,d} + L^\mathrm{intra,max}_{s,\psi(d)} \leq E_s, \quad   &&\forall s\in\mathcal{S},d\in\mathcal{D}, &\\
    & L^\mathrm{inter}_{s,d}(1-\varphi_s)^H + L^\mathrm{intra,min}_{s,\psi(d)} \geq 0, \quad &&\forall s\in\mathcal{S},d\in\mathcal{D}.& 
\end{align}
The number of storage level limit constraints for each storage technology and node changes from $2(D\times H) +P\times H+ D$ in the original formulation (\cref{sec:met_kot_orig}) to $3(P\times H) + 2D$, which results in fewer constraints as long as $P <D\frac{H-1}{H}\approx 350$ representative days.

An almost identical formulation is provided by \citet{Gonzato2021} and \citet{Parolin2024}. In this work, we implement the original formulation by \citet{Kotzur2018} and refer to it as \textit{MinMax}.
\subsection{Storage representation approaches for representative hours}
\label{sec:storage_representation_RH}

The fundamental approach to represent the storage level with RH without relying on chronology is to introduce a new time index for the storage level variable. This new storage time index $j\in\mathcal{J}$ is chronological and thereby rebuilds the sequence of time steps, which was previously lost in the aggregation.

\subsubsection{Full RH storage representation}
\label{sec:full_storage_representation_RH}

\begin{figure}[htbp]
    \centering
    \includegraphics[width=0.9\linewidth]{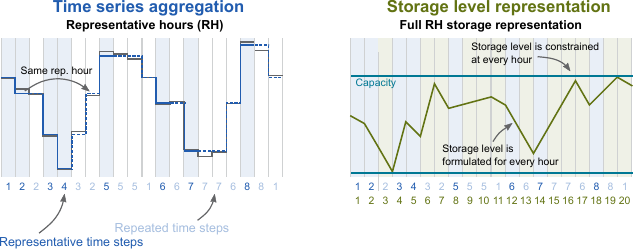}
    \caption{Schematic depiction of the time series aggregation (blue) and storage level representation (green) approaches with representative hours (RH). The fictional fully resolved time series (grey) is aggregated to 8 representative hours. The storage level is represented fully resolved with the method of \cite{Gabrielli2018} (\textit{Full storage resolution}).} 
    \label{fig:schematic_TSA_RH_Gab}
\end{figure}

One option is to fully resolve the storage level variable, i.e., $\vert \mathcal{J}\vert = J = T= 8760$ hours, which is method M1 by \citet{Gabrielli2018} (\cref{fig:schematic_TSA_RH_Gab}). We refer to this method as \textit{Full storage resolution}.

Since $\mathcal{J}=\mathcal{T}$, the original mapping $\sigma : \mathcal{T}\rightarrow\mathcal{I}$ is used to link the storage level and the charging and discharging flows:
\begin{equation}
\label{eq:storage_level_Gab}
    L_{s,j}=L_{s,j-1}(1-\varphi_s)+\Delta H_{s,\sigma(j)}, \quad \forall s\in\mathcal{S},j\in\mathcal{J}\setminus\{j=1\}.
\end{equation}
The other constraints are formulated equivalently:
\begin{align}
\label{eq:add_storage_constraints_Gab}
   & L_{s,j=1}=L_{s,j=J}(1-\varphi_s)+\Delta H_{s,\sigma(j=1)}, \quad && \forall s\in\mathcal{S}, &\\
   & 0\leq L_{s,j} \leq E_s, \quad && \forall s\in\mathcal{S},j\in\mathcal{J}.&
\end{align}
Note that \citet{Gabrielli2018} apply this method to RD; however, due to the generic mapping $\sigma$, it can also be applied to RH.

\subsubsection{Chronological storage level representation}
\label{sec:storage_representation_CRH}
\begin{figure}[htbp]
    \centering
    \includegraphics[width=0.9\linewidth]{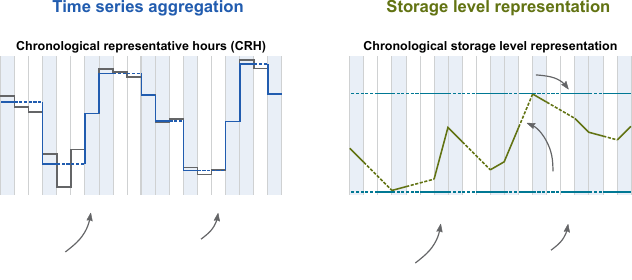}
    \caption{Schematic depiction of the time series aggregation (blue) and storage level representation (green) approaches with chronological representative hours (CRH). The fictional fully resolved time series (grey) is aggregated to 10 representative hours. The storage level is represented by the same number of chronological storage level time steps, as introduced by Pineda et al. \cite{Pineda2018} (\textit{Chrono}). Note that this approach is a special case of the \textit{Proposed method} approach presented in this paper.} 
    \label{fig:schematic_TSA_RH_chron}
\end{figure}

A special case of TSA with RH are chronological representative hours (CRH) \cite{Pineda2018}. Here, we denote this approach as \textit{Chrono}. All time steps that are represented by the same representative hour are adjacent (left plot in \cref{fig:schematic_TSA_RH_chron}). Since the chronology of representative time steps is thereby automatically preserved, the sequence of storage level time steps is equal to the sequence of representative time steps. Note that this relationship between storage time step sequence and representative time step sequence is a special case of the \textit{Proposed method}, where a new storage time step is added every time the representative time step changes. Therefore, the implementation of the \textit{Proposed method} can be used for \textit{Chrono}.







\clearpage
\section{Supplementary Results}

\subsection{Storage level variables and constraints}

\begin{figure}[htbp]
    \centering
    \includegraphics[width=0.8\linewidth]{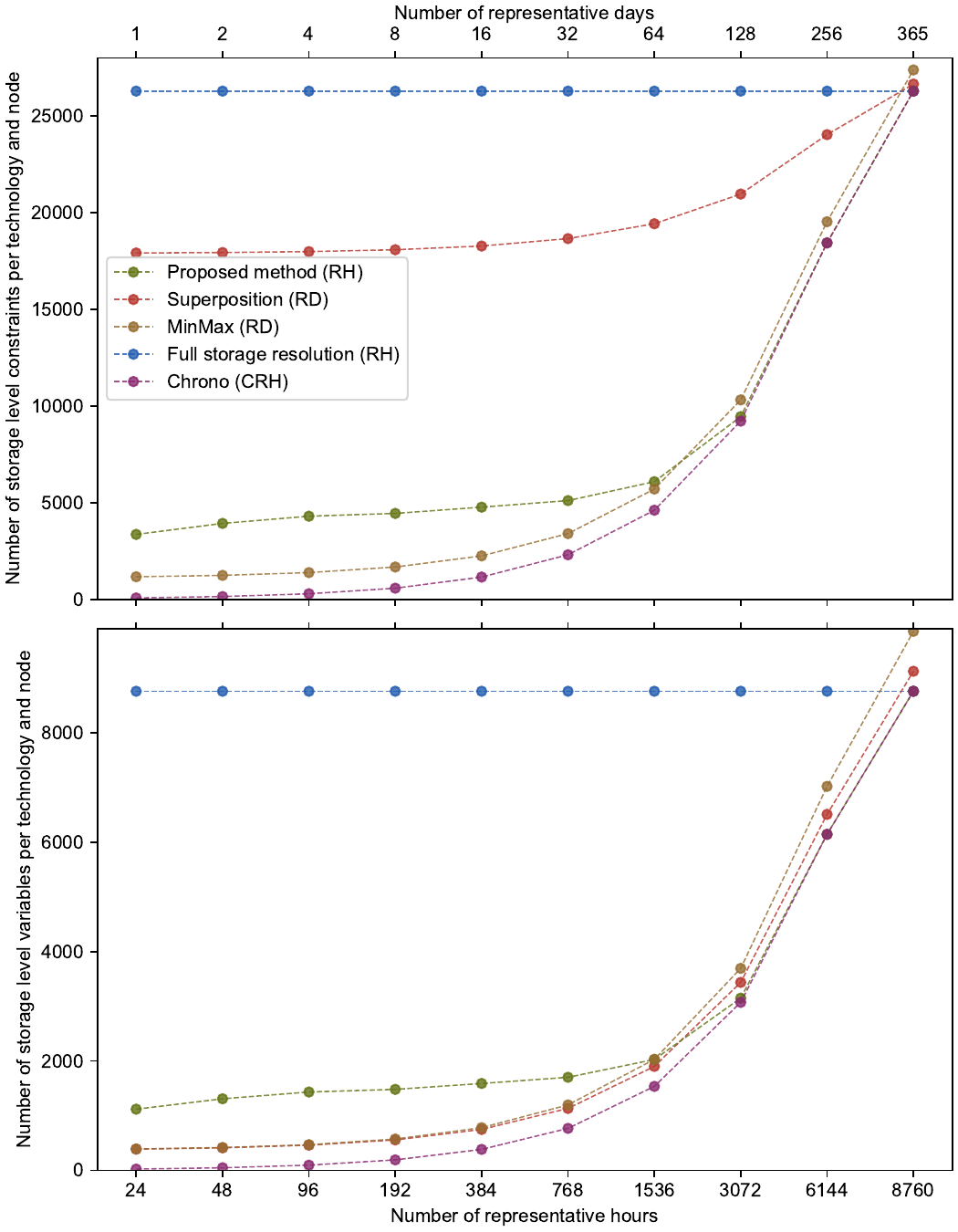}
    \caption{Number of storage level constraints (top) and variables (bottom) for the five investigated storage representation methods. The numbers are reported for each storage technology ($\vert\mathcal{S}\vert=4$) at each node ($\vert\mathcal{N}\vert=28$).} 
    \label{fig:num_cons_vars}
\end{figure}

The differences in solving time can be explained by the size of the optimization problem, in particular, by the number of storage level variables and constraints (\cref{fig:num_cons_vars}). In general, the more storage level variables and constraints are formulated, the larger is the resulting optimization problem, and, therefore, the longer it usually takes to solve \cite{Spielman2003,Wright1997}. 
\textit{Full storage resolution} has a constant number of storage level variables ($T=8760$ per technology and node) and constraints ($3T=26280$ per technology and node). Therefore, \textit{Full storage resolution} shows the highest solving time at the strongest aggregation.

The RD methods \textit{Superposition} and \textit{MinMax} significantly reduce the number of storage level variables compared to \textit{Full storage resolution} to around 400 per technology and node ($365 + 24$ and $365 + 24 + 2$, respectively) for the strongest aggregation of 1 RD. At full resolution, both RD approaches show slightly more storage level variables than \textit{Full storage resolution}. The difference in solving time between the two RD approaches comes from the number of storage level constraints. \textit{Superposition} formulates the limit on storage levels for each original time step (\cref{eq:kotzur_hourly_storage_level_limit}) but does not require a coupling constraint for each original time step; hence, the approach shows fewer constraints than \textit{Full storage resolution} except for the full resolution. \textit{MinMax} requires significantly fewer storage level constraints than \textit{Superposition} for all time step configurations except full resolution, where it requires the most constraints of all approaches.

\begin{figure}[tbph]
    \centering
    \includegraphics[width=0.7\linewidth]{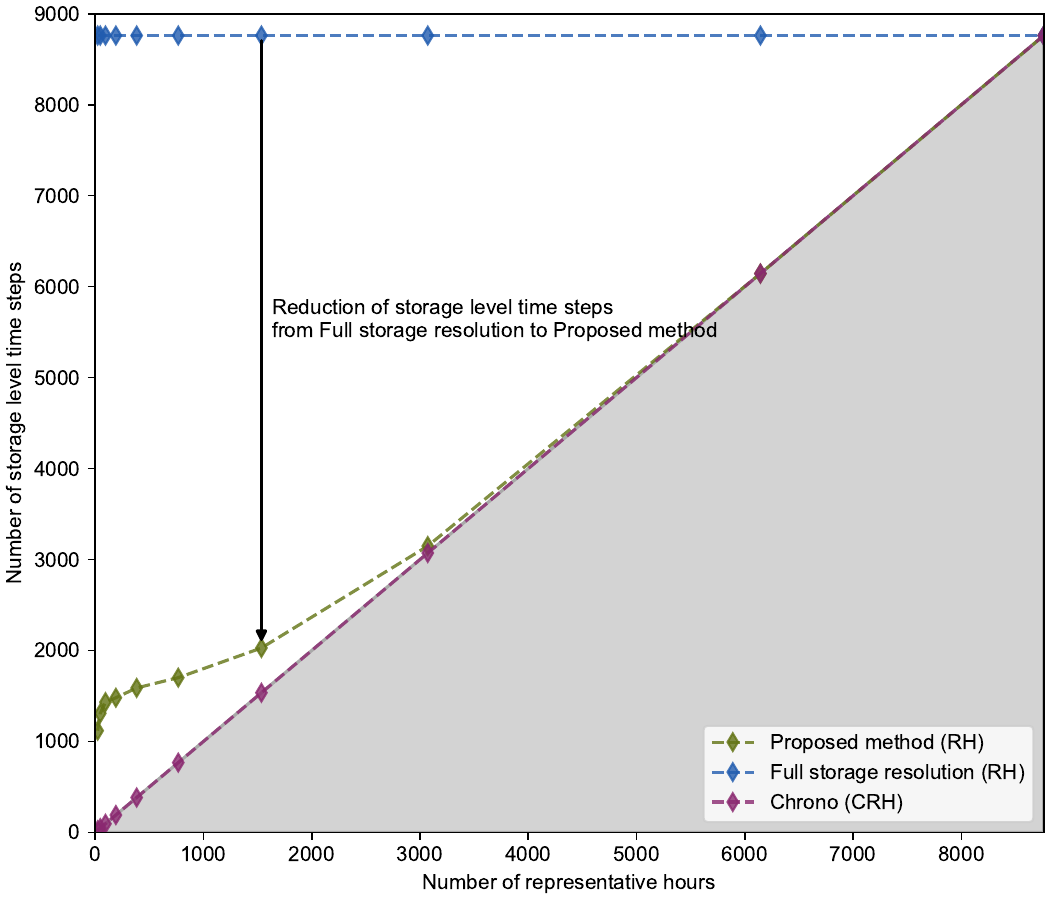}
    \caption{Number of storage level time steps versus number of representative hours for \textit{Full storage resolution}, \textit{Proposed method}, and \textit{Chrono}. Since the mapping $\vartheta:\mathcal{J}\rightarrow\mathcal{I}$ is unique, the number of storage level time steps cannot be lower than the number of representative hours (grey area).} 
    \label{fig:num_stor_ts}
\end{figure}

The remaining RH methods \textit{Proposed method} and \textit{Chrono} show significantly fewer storage level constraints and variables than \textit{Full storage resolution}. \textit{Proposed method} shows more constraints and variables than \textit{MinMax} at strong aggregation, but the numbers do not strongly increase until 1536 RH. Hence, \textit{MinMax} overtakes \textit{Proposed method} in terms of storage level constraints and variables at 3072 RH (128 RD). \textit{Chrono} shows the fewest variables ($I$) and constraints ($3I$), which increase linearly with the number of representative hours. At 3072 RH, \textit{Proposed method} and \textit{Chrono} show a very similar number of constraints and variables.

\cref{fig:num_stor_ts} shows how the RH aggregation results in a different number of storage level time steps. The number of storage level time steps of \textit{Proposed method} is located in the white triangle in \cref{fig:num_stor_ts} between \textit{Full storage resolution} ($J=T=8760$, upper bound) and \textit{Chrono} ($J=I$, lower bound). Remarkably, \textit{Proposed method} is significantly closer to \textit{Chrono} than \textit{Full storage resolution}, without losing any information or accuracy in comparison with \textit{Full storage resolution}. The combinatorial change of representative hours in the sequence $\sigma$ leads to around 1000 storage level time steps at the strongest aggregation of 24 RH. For a higher number of representative hours, the number of storage level time steps increases less strongly to around 2000 at 1536 RH. The number of storage level time steps continues to approach the lower bound of \textit{Chrono}.

\clearpage

\subsection{Electricity price duration curve}
\begin{figure}[htbp]
    \centering
    \includegraphics[width=0.7\linewidth]{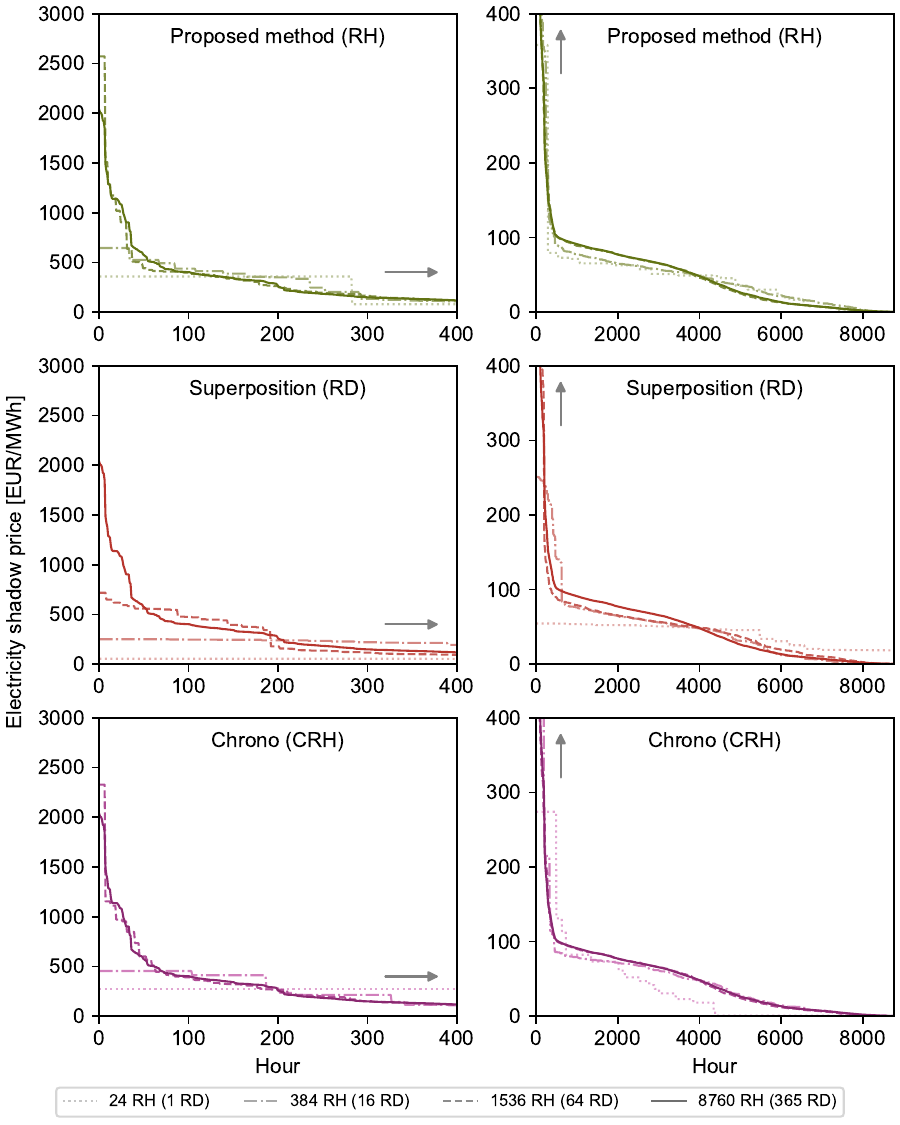}
    \caption{Price duration curve of mean electricity shadow prices (dual variable of energy balance) for \textit{Proposed method} (RH), \textit{Superposition} (RD), and \textit{Chrono} (CRH) for 24, 384, 1536, and 8760 RH. The left column shows the price for the first 400 hours, and the right column shows the price below 400 EUR/MWh. \textit{Full storage resolution} and \textit{MinMax} are not shown, as they have the same system design as \textit{Proposed method} and \textit{Superposition}, respectively.} 
    \label{fig:snapshot_TSA_electricity_price_LDS}
\end{figure}

The solution accuracy is further highlighted by the price duration curve of the electricity shadow prices (\cref{fig:snapshot_TSA_electricity_price_LDS}). The mean electricity price ranges from 0.3 EUR/MWh to 2033 EUR/MWh in the fully resolved model. For 8200 hours (94\% of the year), the electricity price is below 100 EUR/MWh. \textit{Proposed method} and \textit{Chrono} approximate the electricity price well for most TSA settings, whereas \textit{Superposition} shows significantly lower maximum and higher minimum electricity prices. Especially, for the strongest aggregation (24 RH/1 RD), \textit{Superposition} shows a maximum price of 55 EUR/MWh, whereas the maximum price of \textit{Proposed method} is 358 EUR/MWh for the same resolution. \textit{Chrono} shows a similar behavior to \textit{Proposed method} at the highest prices, but the price drops to almost zero for half of the year for 24 RH. Both \textit{Proposed method} and \textit{Chrono} overestimate the maximum electricity prices for 1536 RH.
\clearpage

\bibliographystyle{elsarticle-num-names}

\bibliography{references}
\end{document}